\documentclass{amsart}

\usepackage{amsmath,latexsym,amscd,amsbsy,amssymb,amsfonts,amsthm,fleqn,leqno,
euscript, graphicx, color}

\newtheorem{thm}{Theorem}[section]
\newtheorem{pro}[thm]{Proposition}

\newtheorem{cor}[thm]{Corollary}

\newtheorem{ex}[thm]{Example}
\newtheorem{que}[thm]{Question}

\newtheorem{df}[thm]{Definition}

\newtheorem{rem}[thm]{Remark}

\newcommand{\Int}{\mbox{{\rm int}}\,}
\newcommand{\cl}{\mbox{{\rm cl}}}

\newcommand{\st}{\mbox{{\rm St}}}
\newcommand{\pr}{\mbox{{\rm pr}}}

\newcommand{\id}{\mbox{{\rm id}}}

\begin{document}

\title{\bf  Topological groups that realize homogeneity of topological spaces}
\thanks{The work is supported by NSERC Grant 257231-09}

\author{A.\,V.~Karasev}
\address{Department of Computer Science and Mathematics, Nipissing University,
100 College Drive, P.O. Box 5002, North Bay, ON, P1B 8L7, Canada}
\email{alexandk@nipissingu.ca}

\author{K.\,L.~Kozlov}\thanks{The second author is partially supported by RFBR Grant
15-01-05369}
\address{Department of General Topology and Geometry, Faculty of Mechanics and Mathematics, \\
Moscow State University, Moscow 119991, Russia}
\email{kkozlov@mech.math.msu.su}


\keywords{homogeneity, topological group, $G$-space, uniformity,
completion}

\subjclass{Primary 54H15 Secondary 57S05;\ 22F30;\ 54D35}

\begin{abstract}
We present results on simplifying an acting group preserving
properties of actions: transitivity, being a coset space and
preserving a fixed equiuniformity in case of a $G$-Tychonoff space.
\end{abstract}

\maketitle







\section{Introduction}

A topological space $X$ is {\it homogeneous} if for any $x, y\in X$
there is a homeomorphism $h$ of $X$ such that $h(x)=y$. In the study
of topological homogeneity it is natural to ask about groups (or
their classes) which acts continuously and transitively on a
homogeneous space. The survey of A.\,V.~Arhangel'skii and J.~van
Mill~\cite{ArM} can serve as a good introduction to the subject of
the paper.

Let $G$ acts continuously and transitively on $X$. Using the
conjugation of stabilizes of points from one orbit, a normal
subgroup $K$ of the kernel of action can be chosen. The quotient
group $G/K$ can be equipped with topology (weaker than the quotient
topology) in which the naturally defined action is continuous. This
leads to the possibility to define a continuous and transitive
action on $X$ of a group $H$ such that $\chi(H)\leq\chi(X)\cdot{\rm
inv}(G)$ and ${\rm ib}(H)\leq {\rm ib}(G)$ (Theorem~\ref{th3.2}).
Corollaries~\ref{cor3.10} and~\ref{cor4.22} show what restrictions
we have on a pseudocharacter and an invariance number of a group
which continuously, effectively and transitively acts on a space $X$
of character $\leq\tau$. Let us note that there is a continuous,
effective and transitive action of a metrizable (discrete) group $G$
(of ${\rm ib}(G)\geq |X|$) on a homogeneous space $X$.

We show how to replace transitive actions of $\mathcal G$-range
topological groups to topological groups from class $\mathcal G$
(Theorem~\ref{th3.1}). As a class $\mathcal G$ we can examine groups
of weight or character $\leq\tau$. Transitive actions of
$\tau$-narrow and $\tau$-balanced topological groups on spaces of
character $\leq\tau$ can be replaced to transitive actions of groups
of weight and character $\leq\tau$ respectively
(Corollary~\ref{cor3.5}). In particular, transitive actions of
$\omega$-narrow and $\omega$-balanced topological groups on first
countable spaces can be replaced to transitive actions of separable
metrizable and metrizable groups respectively and transitive action
of a subgroup of a product of \v Cech complete groups can be
replaced to a transitive action of an inframetrizable group
(Corollary~\ref{cor4.1}). As a consequence we have: if an
$\omega$-narrow group acts continuously and transitively on a space
$X$ of countable character with Baire property, then $X$ is a
separable metrizable space (Corollary~\ref{cor4.11}); there is no
continuous and transitive action of an $\omega$-narrow group on the
homogenous space which is not a coset space of Ford
(Example~\ref{ex4.1}), on the two-arrows space, on the homogenous
first countable compactum which is not a coset space of Fedorchuk,
on the Sorgenfrey line (Example~\ref{ex4.3}).

A more strong form of homogeneity is when $X$ is a coset space of a
topological group. Corollaries~\ref{cor3.1} and~\ref{cor4.2} show
that if  $X$ is a coset space of a $\tau$-balanced group (or a
$\tau$-narrow group) and  $\chi (X)\leq\tau$, then $X$ is coset
space of a group of character $\leq\tau$ (or of weight $\leq\tau$).
As a consequence we have: if a subgroup of a product of \v Cech
complete groups acts continuously and transitively on a compactum
$X$ of countable character, then $X$ is a metrizable compactum
(Theorem~\ref{pro4.1}); the Sorgenfrey line  is not a coset space of
an $\omega$-balanced group (Example~\ref{ex4.2}); the two-arrows
space is not a coset space of a subgroup of a product of \v Cech
complete groups (Example~\ref{ex4.4}).

If $X$ is a coset space of a topological group $G$, then for the
natural action $\alpha$ of $G$ on $X$ by left translations the
$G$-space $(G, X, \alpha)$ is $G$-Tychonoff and has
equiuniformities. For continuous transitive actions it is not known.

\medskip

\noindent {\bf Question}~\cite[Question 2.6]{Mg3}. Let $X$ be a
Tychonoff $G$-space with the transitive action. Is it true that $X$
is $G$-Tychonoff?

\medskip

Corollaries~\ref{cor3.21} and~\ref{cor4.21} show that if $X$ is a
coset space of a $\tau$-narrow group and $\chi (X)\leq\tau$, then
there is an equiuniformity on $X$ of weight $\leq\tau$.

The usage of equiuniformities allows to apply the replacing process
preserving properties of actions to non transitive actions. It
consists in equipment of a quotient group with respect to the kernel
of action with the topology of uniform convergence. If $(G, X,
\alpha)$ is a $G$-Tychonoff space and $\mathcal U$ is an
equiuniformity on $X$, then we can replace the group $G$ on an
action of a group $H$ such that $\chi(H)\leq w({\mathcal U})$, ${\rm
ib}(H)\leq {\rm ib}(G)$ and $\mathcal U$ is an equiuniformity
(Theorem~\ref{th3.3}). As a consequence we have: there is no
continuous and effective action on $X$ of a topological group $G$
with $\psi(G)>\tau$ for which $\mathcal U$, $w({\mathcal
U})\leq\tau$, is an equiuniformity (Corollary~\ref{cor3.22}). In
Example~\ref{ex3.7} we make remarks about existence of complete
metrics on a homogeneous Polish space which is not a coset space of
J.~van Mill.


\section{Preliminaries}

All spaces are assumed to be Tychonoff, maps are continuous and
notations, terminology and designations are from~\cite{E}. Nbd(s) is
an abbreviation of open neighborhood(s), $\cl_X A$, $\Int_X A$ are
the closure, interior of a subset $A$ of a space $X$, respectively.
A metrizable space is called Polish if it is separable and has a
complete metric. A map $\id$ is an identity map.

For information about topological groups see~\cite{AT}
and~\cite{RD}. An {\it inframetrizable} group is a subgroup of a \v
Cech complete group~\cite{RD}. By ${\rm ib}(G)$  the {\it index of
narrowness} of a topological group $G$ is denoted and by ${\rm
inv}(G)$ the {\it invariance number} of $G$ is denoted, ${\rm
inv}(G)\leq {\rm ib}(G)$, see~\cite[Chapter 5]{AT}.

By a {\it separating family of homomorphisms} on a topological group
we understand the family of homomorphisms which separates points and
closed sets. By $N_G(e)$ we denote the family of nbds of the unit
$e\in G$. If $\varphi: G\to H$ is an epimorphism then ${\rm ib}
(H)\leq {\rm ib} (G)$. The {\it kernel of homomorphism} $\varphi:
G\to H$ is denoted ${\rm ker}\varphi$.

A class $\mathcal G$ of topological groups is called {\it
$\tau$-multiplicative} if the product of $\tau$ representatives of
$\mathcal G$ is an element of $\mathcal G$. If $\tau=\sigma$ then we
say $\sigma$-multiplicative. If the class is closed under arbitrary
products then it is called {\it multiplicative}.

A class $\mathcal G$ of topological groups is called {\it
hereditary} if a subgroup of a representative of $\mathcal G$ is an
element of $\mathcal G$.


\subsection{$G$-spaces}

By an action of a group $G$ on a set $X$ we mean a map $\alpha:
G\times X\to X$ such that
$$\alpha (e,x)=x,\ \alpha (g,\alpha (h,x))=\alpha (gh,x)$$
for any $x\in X$, $g,h\in G$ (where $e$ is the unit of the group
$G$). We denote by $\alpha_x: G\to X$, $x\in X$, and $\alpha^g: X\to
X$, $g\in G$, the maps defined by the rules $\alpha_x(g)=\alpha
(g,x)=\alpha^g(x)$. If it is clear what action is considered, then
we write $\alpha (g,x):= gx$. The action is {\it transitive} if the
map $\alpha_x$ is a surjection. For a set $X$ with an action
$\alpha$ of a group $G$ the subgroup $H=\{g\in G: gx=x, x\in X\}$ is
called the {\it kernel} of action and is denoted ${\rm ker}_\alpha$.
If ${\rm ker}_\alpha=\{e\}$, then the action is called {\it
effective}.

For a family of subsets $\gamma=\{O_\alpha: \alpha\in A\}$ and $g\in
G$ we set $g\gamma=\{gO_\alpha: \alpha\in A\}$. A family of covers
$\Gamma=\{\gamma_s: s\in S\}$ is {\it saturated} if
$g\gamma_s\in\Gamma$ for any $s\in S$, $g\in G$.

A {\it $G$-space} is a triple  $(G, X, \alpha)$ of a space $X$ with
a fixed continuous action $\alpha: G\times X\to X$ of a topological
group $G$.  The kernel ${\rm ker}_\alpha$ of a continuous action is
a closed subgroup.

Let $(G, X, \alpha_G)$ and $(H, Y, \alpha_H)$ be $G$-spaces. A pair
of maps $(\varphi: G\to H,\ f: X\to Y)$ such that $\varphi: G\to H$
is a homomorphism and the diagram
$$\begin{array}{ccc}
G\times X &\stackrel{\varphi\times f}\longrightarrow & H\times Y \\
\downarrow\lefteqn{\alpha_G} && \downarrow\lefteqn{\alpha_H} \\
X &\stackrel{f}\longrightarrow & Y
\end{array}$$
is commutative is called an {\it equivariant pair of maps of
$G$-spaces}. Designation $(\varphi, f): (G, X, \alpha_G)\to (H, Y,
\alpha_H)$. The commutativity of the diagram may be written as the
fulfillment of the following condition
$$f(gx)=\varphi(g)f(x)\ \mbox{for any}\ x\in X,\ g\in G.$$
If $f$ is an embedding then the  equivariant pair $(\varphi, f)$ is
an {\it equivariant embedding} of $(X, G_X, \alpha_X)$ into $(Y,
G_Y, \alpha_Y)$. The composition of equivariant pairs of maps is an
equivariant pair.

A $G$-space $(G, X, \alpha)$ is called {\it $G$-Tychonoff}, if there
is an equivariant embedding $(\id, f)$ of $(G, X, \alpha)$ into a
$G$-space $(G, bX, \tilde\alpha)$ where $bX$ is a compactification
of $X$ ($(G, bX, \tilde\alpha)$ (together with the embedding $(\id,
f)$) is a {\it $G$-compactification} of $(G, X, \alpha)$).  The {\it
maximal $G$-compactification} is denoted $(G, \beta_GX,
\alpha_\beta)$.


\subsection{$d$-Open actions}

Information about ($d$-)open actions, their properties and natural
uniform structures which are generated by ($d$-)open actions can be
found in~\cite{CK3}, \cite{K4} and \cite{K2013}. For the convenience
of the reader we remind some facts which are used below.

\begin{df}{\rm (see~\cite{CK3})}\label{d2.2.1}
{\rm An action  $\alpha: G\times X\to X$ is called
\begin{itemize}
\item[] {\it open} if $x\in\Int (Ox)$ for any point $x\in X$ and  $O\in N_G(e)$;
\item[] {\it $d$-open} if $x\in\Int (\cl (Ox))$ for any point $x\in X$ and  $O\in
N_G(e)$.
\end{itemize}}
\end{df}

\begin{rem}\label{r2.2.1}
{\rm Terminology in definition~\ref{d2.2.1} is motivated by the fact
that an action is ``open" (``$d$-open") iff maps $\alpha_x:G\to X$,
$x\in X$, are open ($d$-open)~\cite[Remark 1.7]{K2013}.

F.~Ancel~\cite{A} called continuous open actions {\it
micro-transitive} and continuous $d$-open actions {\it weakly
micro-transitive}.

If $(G, X, \alpha)$ is with a $d$-open action then $X$ is a direct
sum of clopen subsets ({\it components of the action}). Each
component of the action is the closure of the orbit of an arbitrary
point of this component. If an action is open then $X$ is a direct
sum of clopen subsets which are the orbits of the
action~\cite[Remark 1.7]{K2013}. For a $G$-space $X$ with an open
action and one component of action $X$ is the quotient space of
$G$.}
\end{rem}

The following theorem is a generalization of a correspondent result
of F.~Ancel for complete metrizable groups~\cite{A}.

\begin{thm}{\rm \cite[Theorem 3]{K4}}\label{th2.2.1}
Let $(G, X, \alpha)$ be a $G$-space with a $d$-open action and let
for any $x\in X$ the coset space $G/G_x$ be \v Cech complete. Then
the action is open.

In particular, a $d$-open continuous action of a \v Cech complete
group is open.
\end{thm}


\subsection{Uniform structures on $G$-spaces}

Uniform structures on spaces are introduced by the families of
covers~\cite{Is} and are compatible with their topology.  For covers
$u$ and $v$ of $X$ we denote $u\succ v$, $u\ast\succ v$ if $u$
refines, respectively star-refines $v$.

A uniformity $\mathcal U$ on a $G$-space $X$ is called {\it
equiuniformity} (see, for example, \cite{Mg}) if it is saturated and
{\it bounded} (i.e. that for any $u\in\mathcal U$ there exist $O\in
N_G(e)$ and $v\in\mathcal U$ such that $Ov=\{OV: V\in v\}\succ u$).
The action on a $G$-Tychonoff space $(G, X, \alpha)$ with
equiuniformity $\mathcal U$ can be extended to the continuous action
of the Raikov completion $\rho G$ of $G$ on the completion $\tilde
X$ of $X$ with respect to the equiuniformity $\mathcal
U$~\cite[Theorem 3.1]{Mg}. By {\it equivariant completion} of  $(G,
X, \alpha)$ we understand a $G$-space $(G, \tilde X,
\tilde\alpha_\rho)$ with the natural equivariant embedding $(\id,
i): (G, X, \alpha)\to (G, \tilde X, \tilde\alpha)$.

\begin{thm}{\rm \cite{CK1}}\label{th2.3.1''}
If the action on a $G$-space $(G, X, \alpha)$ is $d$-open, then the
family of covers $\gamma_O=\{\Int (\cl (Ox)): x\in X\}$, $O\in
N_G(e)$, is the base of the maximal equiuniformity on $X$.
\end{thm}

\subsection{Topologizing groups using equiniformities on $G$-Tychonoff spaces}

\begin{thm}\cite[Theorem 2.2]{Ford}\label{thm6.1}
If $X$ possesses a uniform structure $\mathcal U$ under which every
element of some homeomorphism group $G$ is uniformly continuous,
then $G$ is a topological group relative to the uniform convergence
notion induced by the uniformity.
\end{thm}

A bijection $f$ of uniform spaces is called a {\it uniform
equivalence} if $f$ and $f^{-1}$  are uniformly continuous.

\begin{rem}\label{r6.1}
{\rm The sets $O_u=\{g\in G: \forall x\in X\ g(x)\in\st(x, u)\}$,
$u\in \mathcal U$, are a base of nbds (not open) at $e$.

The uniformity of uniform convergence in theorem~\ref{thm6.1}
coincides with the right uniformity on $G$~\cite[Chapter 2, Exercise
2]{RD}\label{rA.1}.}
\end{rem}

\begin{pro}\cite[Remark 4.2 (c)]{AAK}\label{pro6.1}
Uniformity $\mathcal U$ in theorem~\ref{thm6.1} is an equiuniformity
and, hence, $(G, X, \alpha)$ is a $G$-Tychonoff space.
\end{pro}

\begin{cor}\cite[Corollary 4.5]{AAK}\label{c6.1}
If $w(\mathcal U)\leq\aleph_0$ for a totally bounded uniformity
$\mathcal U$ on $X$ then the group of uniform equivalences $H(X,
{\mathcal U})$ in the topology of uniform convergence is separable
and metrizable.
\end{cor}

The following theorem is a generalization of~\cite[theorem
4.6]{AAK}.

\begin{thm}\cite[Lemma 5]{Mart}\label{thm6.2}
Let $(G, X, \alpha)$ be a $G$-Tychonoff space with an effective
action, $\mathcal U$ be an equiuniformity on $X$. Then
\begin{enumerate}
\item each element
of $G$ is a uniform equivalence (with respect to $\mathcal U$),
\item the topology of uniform convergence $\tau_{\mathcal U}$ on $G$ is
coarser than the original one,
\item $((G, \tau_{\mathcal U}), X, \alpha)$
is a $G$-Tychonoff space and $\mathcal U$ is an equiuniformity on
$X$.
\end{enumerate}
\end{thm}

This result shows, in fact, the following.

\begin{cor}
Let $\mathcal U$ be an uniformity on $X$. For any topological group
$G$ such that there is an effective action $G\times X\to X$ for
which $(G, X, \alpha)$ is a $G$-Tychonoff space and $\mathcal U$ is
an equiuniformity on $X$, the topology of uniform convergence on $G$
is coarser than the original topology on $G$.
\end{cor}


\subsection{Properties of equivariant pair of maps}

\begin{pro}\label{pro2.1}
Let $(\varphi, {\rm id}): (G, X, \alpha)\to (H, X, \gamma)$ be an
equivariant pair of maps of $G$-spaces and $\varphi: G\to H$ be an
epimorphism. Then we have
\begin{itemize}
\item[(a)] ${\rm ker}\ \varphi\subset{\rm ker}_\alpha$;
\item[(b)] $\varphi({\rm ker}_\alpha)={\rm ker}_\gamma$ and if $\alpha$ is effective then $\gamma$ is effective;
\item[(c)] $Gx=Hx$ for any $x\in X$, in particular, if $\alpha$ is transitive then $\gamma$ is transitive;
\item[(d)] if $\alpha$ is open ($d$-open) then $\gamma$ is open ($d$-open) and components of actions $\alpha$ and $\gamma$ coincide;
\item[(e)] if $\mathcal U$ is an equiuniformity on $X$ in $(H, X, \gamma)$ then $\mathcal U$ is an equiuniformity on $X$ in $(G, X,
\alpha)$, hence, if $(H, X, \gamma)$ is $G$-Tychonoff then $(G, X,
\alpha)$ is $G$-Tychonoff.
\end{itemize}
\end{pro}

{\bf Proof.}\quad Properties (a), (b) and (c) immediately follows
from the equality $\alpha(g, x)=\gamma(\varphi(g), x)$, $x\in X$,
$g\in G$.

In order to prove property (d) take $x\in X$ and $O\in N_H(e)$. For
$U=\varphi^{-1}(O)\in N_G(e)$
$$x\in\Int(Ux)=\Int(\varphi(U)x)=\Int(Ox)\ \mbox{in case of an open action}\ \alpha,$$
$$x\in\Int(\cl(Ux))=\Int(\cl(\varphi(U)x))=\Int(\cl(Ox))\ \mbox{in case of a}\ d\mbox{-open action}\ \alpha.$$
Thus openness ($d$-openness) is preserved. Since
$\Int(\cl(Gx))=\Int(\cl(Hx))$, $x\in X$, the components of actions
$\alpha$ and $\gamma$ coincide.

Property (e). If $\mathcal U$ is an equiuniformity on $X$ in $(H, X,
\gamma)$ then it is a uniformity on $X$. From the equality
$\alpha(g, A)=\gamma(\varphi(g), A)$ for any $A\subset X$, $g\in G$
it follows that $\mathcal U$ is a saturated uniformity in $(G, X,
\alpha)$. Its boundedness follows from the equality
$\alpha(\varphi^{-1}(O), x)=\gamma(O, x)$, $x\in X$, $O\in N_H(e)$.

The last statement in property (e) is a straightforward
consequence.\ $\Box$

\begin{rem}
{\rm If $G$ is a discrete group then a $G$-space $(G, X, \alpha)$ is
$G$-Tychonoff. This yields the existence of an example of a
$G$-Tychonoff space which image under an equivariant pair of maps is
not $G$-Tychonoff.}
\end{rem}

From proposition~\ref{pro2.1} (e) we have.

\begin{cor}\label{cor2.1}
Let $(G, X, \alpha)$ be a $G$-Tychonoff space and $\mathcal U$ be
the maximal equiuniformity on $X$. If $(\varphi, {\rm id}): (G, X,
\alpha)\to (H, X, \gamma)$ is an equivariant pair of maps and
$\mathcal U$ is an equiuniformity on $X$ in $(H, X, \gamma)$, then
$\mathcal U$ is the maximal equiuniformity on $X$ in $(H, X,
\gamma)$.
\end{cor}


\section{Replacing acting group (general case)}

Let $\mathcal G$ be a class of topological groups. A topological
group $G$ is called {\it range-$\mathcal G$} if for any $O\in
N_G(e)$ there exists a continuous homomorphism $h$ of $G$ to
$H\in{\mathcal G}$ such that $h^{-1}(U)\subset O$ for some $U\in
N_H(e)$~\cite[\S\ 3.4]{AT}.

\begin{rem}\label{rem3.1}~\cite[Theorem 3.4.21]{AT} {\rm Let $\mathcal G$ be a class of topological groups. The following conditions for a group
$G$ are equivalent:
\begin{itemize}
\item $G$ is range-$\mathcal G$;
\item $G$ has a separating
family of homomorphisms to topological groups from $\mathcal G$;
\item $G$ is topologically isomorphic to a subgroup of
a product of a family of topological groups from $\mathcal G$.
\end{itemize}}
\end{rem}

\begin{pro}\label{pro3.1}
Assume that $(G, X, \alpha)$ is a $G$-space and a group $G$ is a
subgroup of the product $\Pi=\Pi\{G_s\in{\mathcal G}: s\in S\}$ of a
family of topological groups. If

\noindent a) the action $\alpha$ is transitive and $\chi
(X)\leq\tau$ or

\noindent b)  $(G, X, \alpha)$ is $G$-Tychonoff, $\mathcal U$ is an
equiuniformity on $X$ and $w({\mathcal U})\leq\tau$,

\noindent then there exist: $S'\subset S$ such that $|S'|\leq\tau$;

\noindent an action $\gamma:\pr_{S'}(G)\times X\to X$, where
$\pr_{S'}:\Pi\to\Pi_{S'}=\{G_s\in{\mathcal G}: s\in S'\}$ is a
projection, such that $(\pr_{S'}(G), X, \gamma)$ is a $G$-space with
transitive action $\gamma$ in case a) or $\mathcal U$ is an
equiuniformity on $X$ in $(\pr_{S'}(G), X, \gamma)$ in case b); and

\noindent an equivariant pair of maps $(\varphi=\pr_{S'}|_G, {\rm
id}):(G, X, \alpha)\to (\pr_{S'}(G), X, \gamma)$.
\end{pro}

{\bf Proof.}\quad In case a) fix a point $x\in X$ and its base of
nbds $\{U_t: t\in T\}$, $|T|\leq\tau$. For each $U_t$ take $O_t\in
N_G(e)$ such that $O_tU_r\subset U_t$ for some nbd $U_r$ of $x$.

In case b) fix a base $\{U_t: t\in T\}$, $|T|\leq\tau$, of $\mathcal
U$. For each $U_t$ take $O_t\in N_G(e)$ such that $\{O_tx: x\in
X\}\succ U_t$.

In both cases, without loss of generality, we may examine each $O_t$
as the trace on $G$ of a rectangular set in $\Pi$ which depends on
the finite number of coordinates $S_t$, $t\in T$. The cardinality of
the set $S'=\bigcup\{S_t: t\in T\}$ is $\tau$.

Let $\Pi_{S'}=\Pi\{G_s\in{\mathcal G}: s\in S'\}$,
$\pr_{S'}:\Pi\to\Pi_{S'}$, $H=\pr_{S'}(G)\subset\Pi_{S'}$.
Evidently, $\varphi=\pr_{S'}|_G$ is a continuous homomorphism of $G$
onto $H$, and $\varphi (O_t)$ is open in $H$,  $\varphi^{-1}(\varphi
(O_t))=O_t$ for all $t\in T$.

The straightforward checking allows to prove the following
properties:
\begin{itemize}
\item $\mbox{for any}\ O_t\ \mbox{there is}\ U(t)\in N_H(e)\
\mbox{such that}\ \varphi^{-1}(U(t))\subset O_t, t\in T,$
\item $\mbox{for any}\ U\in N_H(e)\ \mbox{and}\ g\in G\ \mbox{there is}\
V\in N_H(e)\ \mbox{such that}\
g^{-1}\varphi^{-1}(V)g\subset\varphi^{-1}(U).$
\end{itemize}

These properties yields that the action $\gamma$ of $H$ on $X$ (for
$y\in X$, $g\in H$ take $h\in\varphi^{-1}(g)$ and put $\gamma (g,
y)=\alpha (h, y)$) is well defined. In fact, by the definition of
$H$ for any $h\in\ker\varphi$ we have $\alpha (h, x)=x$. For $y=gx$
and $h\in\ker\varphi$ we have $$\alpha (h, gx)=\alpha (hg, x)=\alpha
(gh', x)=\alpha (g, h'x)=\alpha (g, x)=gx\ \mbox{for}\
h'=g^{-1}hg\in\ker\varphi.$$

Further, $(H, X, \gamma)$ is a $G$-space in case a). In fact, every
$\gamma_h: X\to X$, $h\in H$, is continuous. For $y=gx$ and its nbd
$W$ the set $g^{-1}W$ is a nbd of $x$. There is a nbd $U$ of $x$ and
$O=\varphi^{-1}(\varphi (O))$  such that $OU\subset g^{-1}W$. The
set $gU$ is a nbd of $gx$. There is $V=\varphi^{-1}(\varphi (V))$
such that $g^{-1}Vg\subset O$. Then $VgU\subset
gOg^{-1}gU=gOU\subset W$.

In case b) $(H, X, \gamma)$ is a $G$-Tychonoff space and $\mathcal
U$ is an equiuniformity on $X$ in  $(H, X, \gamma)$. In fact, every
$\gamma_h: X\to X$, $h\in H$, is uniformly continuous with respect
to $\mathcal U$ and for every $u, v\in\mathcal U$, $v\star\star\succ
u$, there is $O\in N_H(e)$ such that $\{\varphi^{-1}(O)x=Ox: x\in
X\}\succ v$ and $Ov\succ u$. This yields that the action $\gamma$ is
continuous and bounded by $\mathcal U$.

Evidently, $(\varphi, {\rm id}):(G, X, \alpha)\to (H, X, \gamma)$ is
an equivariant pair of maps in both cases. Therefore the action
$\gamma$ is transitive in case a)~\ref{pro2.1} (c).\ $\Box$

\medskip

From proposition~\ref{pro3.1} and remark~\ref{rem3.1} we have.

\begin{thm}\label{th3.1}
Assume that $\mathcal G$ is a $\tau$-multiplicative and hereditary
class of topological groups and $G$ is range-$\mathcal G$. Let $(G,
X, \alpha)$ be

\noindent a) a $G$-space with a transitive action and $\chi
(X)\leq\tau$ or

\noindent b) a $G$-Tychonoff space with an equiuniformity $\mathcal
U$ on $X$ and $w({\mathcal U})\leq\tau$.

Then there exist: a $G$-space $(H, X, \gamma)$ with a transitive
action $ \gamma$ in case a) or a $G$-Tychonoff space $(H, X,
\gamma)$ with the equiuniformity $\mathcal U$ such that
$H\in\mathcal G$; and

\noindent an equivariant pair of maps $(\varphi, {\rm id}):(G, X,
\alpha)\to (H, X, \gamma)$, where $\varphi$ is an epimorphism.
\end{thm}

\begin{rem}\label{rem3.2}{\rm
The following $\tau$-multiplicative and hereditary classes $\mathcal
G$ of topological groups and corresponding range-$\mathcal G$
classes of groups are well-known~\cite[\S\ 5.1]{AT}.

\medskip

\begin{tabular}{l|llllll}
class $\mathcal G$ of groups & range-$\mathcal G$ class of
groups \\
\hline
groups of weight $\leq\tau$  &  $\tau$-narrow groups \\
& (or groups $G$ with ${\rm ib}(G)\leq\tau$)  \\
 & \\
groups of character $\leq\tau$ & $\tau$-balanced groups (or groups $G$ with ${\rm inv}(G)\leq\tau$)  \\
\end{tabular}

\medskip

The classes of $\tau$-narrow and $\tau$-balanced group are
multiplicative  and hereditary.

Every $\tau$-narrow group is $\tau$-balanced.}
\end{rem}

From property (c) of proposition~\ref{pro2.1}, theorem~\ref{th3.1}
and remark~\ref{rem3.2} we have.

\begin{cor}\label{cor3.5}
Let $(G, X, \alpha)$ be a $G$-space, $\alpha$ be a transitive
action, $\chi(X)\leq\tau$. If $G$ is a $\tau$-narrow (respectively
$\tau$-balanced) group, then there exists a topological group $H$
such that $w(H)\leq\tau$ (respectively $\chi(H)\leq\tau$) and $H$
admits a continuous transitive action on $X$.
\end{cor}

From theorem~\ref{th3.1}, remark~\ref{rem3.2} and
corollary~\ref{cor2.1} we have.

\begin{cor}\label{cor3.6}
Let $(G, X, \alpha)$ be a $G$-Tychonoff space and $\mathcal U$ be an
equiuniformity (the maximal equiuniformity) on $X$ such that
$w({\mathcal U})\leq\tau$. Then there exist: a $G$-Tychonoff space
$(H, X, \gamma)$ such that
\begin{itemize}
\item[(a)]  $w(H)\leq\tau$ if $G$ is a $\tau$-narrow group,
\item[(b)]  $\chi(H)\leq\tau$ if $G$ is a $\tau$-balanced group,
\end{itemize}
$\mathcal U$ is an equiuniformity (the maximal equiuniformity) on
$X$ in $(H, X, \gamma)$; and

\noindent an equivariant pair of maps $(\varphi, {\rm id}):(G, X,
\alpha)\to (H, X, \gamma)$, where $\varphi$ is an epimorphism.
\end{cor}

From property (d) of proposition~\ref{pro2.1}, theorem~\ref{th3.1}
and remark~\ref{rem3.2} we have.

\begin{cor}\label{cor3.1}
Let $\chi(X)\leq\tau$.
\begin{itemize}
\item[(a)] If $X$ is a coset space of a $\tau$-balanced group (or there is a transitive and $d$-open action on $X$ of a
$\tau$-balanced group), then $X$ is a coset space of a topological
group $H$ (or there exists a transitive and $d$-open action of a
group $H$) such that $\chi(H)\leq\tau$.
\item[(b)] If $X$ is a coset space of a $\tau$-narrow group (or there is a transitive and $d$-open action on $X$ of a
$\tau$-narrow group), then $X$ is a coset space of a topological
group $H$ (or there exists a transitive and $d$-open action of a
group $H$) such that $w(H)\leq\tau$ and, hence, $w(X)\leq\tau$.
\end{itemize}
Moreover, there exists an equivariant pair of maps $(\varphi, {\rm
id}):(G, X, \alpha)\to (H, X, \gamma)$, where $\varphi$ is an
epimorphism.
\end{cor}

\begin{cor}\label{cor3.21}
Let $\chi(X)\leq\tau$. If $X$ is a coset space of a $\tau$-balanced
group $G$ (or there is a transitive and $d$-open action on $X$ of a
$\tau$-balanced group $G$), then there exists an equiuniformity
$\mathcal U$ on $X$ in  $(G, X, \alpha)$ such that $w({\mathcal
U})\leq\tau$.

Let $w(X)\leq\tau$. If $X$ is a coset space of a $\tau$-balanced
group $G$ (or there is a transitive and $d$-open action on $X$ of a
$\tau$-balanced group $G$), then there exist:

\noindent an equivariant extension $(G, \tilde X, \tilde\alpha)$ of
$(G, X, \alpha)$ such that $w(\tilde X)\leq\tau$; and

\noindent a complete equiuniformity $\tilde{\mathcal U}$ on $\tilde
X$ such that $w(\tilde{\mathcal U})\leq\tau$.
\end{cor}

{\bf Proof.}\quad By corollary~\ref{cor3.1} there are: an open (or
$d$-open) action on $X$ of a group $H$ such that $\chi(H)\leq\tau$
and

\noindent an equivariant pair of maps $(\varphi, {\rm id}):(G, X,
\alpha)\to (H, X, \gamma)$, where $\varphi$ is an epimorphism.
By~\cite[Proposition 4]{CK3} the weight of the maximal
equiuniformity $\mathcal U$ on $X$ in $(H, X, \gamma)$ is
$\leq\tau$. By property (e) of proposition~\ref{pro2.1} $\mathcal U$
is an equiuniformity on $X$ in $(G, X, \alpha)$.

The proof of the second statement. Since $\chi(X)\leq w(X)\leq\tau$
there exists an equiuniformity $\mathcal U$ on $X$ in  $(G, X,
\alpha)$ such that $w({\mathcal U})\leq\tau$. The completion of $X$
with respect to $\mathcal U$ is the required equivariant completion
$\tilde X$ of $X$ of weight $\tau$ and by~\cite[Theorem 3.1]{Mg}
there exists a natural extension $\tilde\alpha: G\times\tilde
X\to\tilde X$ of $\alpha$.\ $\Box$

\begin{thm}\label{th3.2}
Let $(G, X, \alpha)$ be a $G$-space, $\alpha$ be a transitive
action.

Then there exist a $G$-space $(H, X, \gamma)$ such that $\gamma$ is
a transitive action,
$$\chi(H)\leq\chi(X)\cdot{\rm inv}(G),$$
$${\rm ib}(H)\leq {\rm ib}(G),$$
$$w(H)\leq\chi(X)\cdot{\rm ib}(G),$$
and an equivariant pair of maps $(\varphi, {\rm id}):(G, X,
\alpha)\to (H, X, \gamma)$, where $\varphi$ is an epimorphism.
\end{thm}

{\bf Proof.}\quad The group $G$ is a subgroup of the product
$\Pi=\Pi\{G_s\in{\mathcal G}: s\in S\}$ of a family of topological
groups of character $\leq {\rm inv}(G)$~\cite[Theorem 5.1.9]{AT}. By
proposition~\ref{pro3.1} there exist $S'\subset S$ such that
$|S'|\leq\chi (X)$,

\noindent $(H=\pr_{S'}(G), X, \gamma)$ is a $G$-space with a
transitive action $\gamma$ and

\noindent an equivariant pair of maps $(\varphi, {\rm id}):(G, X,
\alpha)\to (H, X, \gamma)$, where $\varphi=\pr_{S'}|_G$.

It is easy to see that $\chi (H)\leq\chi(X)\cdot{\rm inv}(G)$ and
${\rm ib}(H)\leq {\rm ib}(G)$. Hence,
$$w(H)=\chi(H)\cdot{\rm ib}(H)\leq\chi(X)\cdot{\rm inv}(G)\cdot{\rm
ib}(G)=\chi(X)\cdot{\rm ib}(G).\ \Box$$

\begin{cor}\label{cor3.10}
Let $\chi(X)\leq\tau$. There is no continuous, effective and
transitive action on $X$ of a topological group $G$ with
$\psi(G)>\tau$, ${\rm inv}(G)\leq\tau$.
\end{cor}

{\bf Proof.}\quad Let $G$ acts continuously, effectively and
transitively on $X$ and ${\rm inv}(G)\leq\tau$. Then by
theorem~\ref{th3.2} there exists a group $H$ such that $\chi
(H)\leq\tau$, $H$ acts continuously, effectively and transitively on
$X$ and $H$ is a continuous homomorphic image of $G$. By property
(a) of proposition~\ref{pro2.1} the map of $G$ onto $H$ is
one-to-one. Then, $\psi(G)\leq\tau$.\ $\Box$

\begin{thm}\label{th3.3}
Let $(G, X, \alpha)$ be a $G$-Tychonoff space, $\mathcal U$ be an
equiuniformity on $X$.

Then there exist a $G$-Tychonoff space $(H, X, \gamma)$ such that
$\mathcal U$ is an equiuniformity on $X$,
$$\chi(H)\leq w({\mathcal U}),$$
$${\rm ib}(H)\leq {\rm ib}(G),$$
$$w(H)\leq w({\mathcal U})\cdot{\rm ib}(G),$$
and an equivariant pair of maps $(\varphi, {\rm id}):(G, X,
\alpha)\to (H, X, \gamma)$, where $\varphi$ is an epimorphism.
\end{thm}

{\bf Proof.}\quad Let us assume that the action $\alpha$ is
effective. Otherwise, take a $G$-Tychonoff space $(G/{\rm
Ker}_\alpha, X, \alpha')$.

From theorems\ref{thm6.2} the topology of uniform convergence
$\tau_{\mathcal U}$ on $G$ is weaker than the original topology,
$H=(G, \tau_{\mathcal U})$ is a topological group by
theorem~\ref{thm6.1}, and the naturally defined action $H\times X\to
X$ is continuous. Therefore, $(H, X, \gamma)$ is a $G$-Tuchonoff
space, $\mathcal U$ is an equiuniformity on $X$ and the equivariant
pair of maps of $(G, X, \alpha)$ to $(H, X, \gamma)$ is well
defined. Thus, ${\rm ib}(H)\leq {\rm ib}(G)$.

Since the sets $O_u=\{g\in G: \forall x\in X\ g(x)\in\st(x, u)\}$,
$u\in \mathcal U$, are a base of nbds (not open) at $e$,
$\chi(H)\leq w({\mathcal U})$. The last inequality follows
from~\cite[Theorem 5.2.3]{AT}.\ $\Box$

\begin{rem}{\rm Let us note that if the action $\alpha: G\times X\to
X$ in theorem~\ref{th3.3} is transitive, (d-)open, then the action
$\gamma: H\times X\to X$ is transitive, (d-)open by
proposition~\ref{pro2.1}.}
\end{rem}

\begin{cor}\label{cor3.22}
Let $\mathcal U$ be a uniformity on $X$ such that $w({\mathcal
U})\leq\tau$. There is no continuous and effective action on $X$ of
a topological group $G$ with $\psi(G)>\tau$ for which $\mathcal U$
is an equiuniformity on $X$.
\end{cor}

{\bf Proof.}\quad Let $G$ acts continuously and effectively on $X$,
$\mathcal U$ be an equiuniformity on $X$ such that $w({\mathcal
U})\leq\tau$. Then by theorem~\ref{th3.3} there exists a group $H$
such that $\chi (H)\leq\tau$, $H$ acts continuously and effectively
on $X$,  $(H, X, \gamma)$ is a $G$-Tychonoff space, $\mathcal U$ is
an equiuniformity on $X$ and $H$ is a continuous homomorphic image
of $G$. By property (a) of proposition~\ref{pro2.1} the map of $G$
onto $H$ is one-to-one. Then $\psi(G)\leq\tau$.\ $\Box$


\section{Replacing acting group (countable case)}

\begin{rem}\label{rem3.3}{\rm
The following $\sigma$-multiplicative and hereditary classes
$\mathcal G$ of topological groups and corresponding range-$\mathcal
G$ classes of topological groups are well-known~\cite[\S\ 3.4]{AT},
\cite[Chapter 13]{RD}.

\medskip

\begin{tabular}{l|llllll}
class $\mathcal G$ of groups & range-$\mathcal G$ class of
groups \\
\hline
groups of countable weight  & $\omega$-narrow groups  \\
(or separable metrizable groups) &   \\
groups of countable character  & $\omega$-balanced groups  \\
(or metrizable groups) &  \\
inframetrizable groups & subgroups of products of \v Cech complete groups  \\
(or subgroups of \v Cech complete groups) &  \\
\end{tabular}

\medskip

The classes of  $\omega$-narrow groups, $\omega$-balanced groups and
subgroups of products of \v Cech complete groups are
$\sigma$-multiplicative  and hereditary.

Every $\omega$-narrow group is $\omega$-balanced. Every
$\omega$-balanced group is a subgroup of a product of \v Cech
complete groups.}
\end{rem}

From condition (c) of proposition~\ref{pro2.1}, theorem~\ref{th3.1}
and remark~\ref{rem3.3} we have.

\begin{cor}\label{cor4.1}
Let $(G, X, \alpha)$ be a $G$-space, $\alpha$ be a transitive
action, $\chi(X)\leq\aleph_0$.  If $G$ is an $\omega$-narrow group
(respectively $\omega$-balanced  group or a subgroup of a product of
\v Cech complete groups), then there exists a separable metrizable
(respectively metrizable or inframetrizable) group $H$ such that $H$
admits a continuous transitive action on $X$.

Moreover, in the case of a transitive action of an $\omega$-narrow
group, $X$ is a separable space with countable cellularity.
\end{cor}

\begin{ex}\label{ex4.1}
{\rm It is easy to show that the cellularity of a homogeneous space
$X$ which is not a coset space from~\cite[\S\ 5 Two examples]{Ford}
is uncountable and $\chi (X)=\aleph_0$. If $X$ admits a transitive
continuous action of an $\omega$-narrow group then by
corollary~\ref{cor4.1} $X$ admits a transitive continuous action of
a separable metrizable group and must have a countable cellularity.
We obtain a contradiction. Hence, there is no continuous and
transitive action of an $\omega$-narrow group on $X$.}
\end{ex}

From condition (d) of proposition~\ref{pro2.1}, theorem~\ref{th3.1},
remark~\ref{rem3.3} and~\cite[Corollary 4]{CK3} (which, among other
things, claimes metrizability of a phase space $X$ of a $G$-space
$(G, X, \alpha)$, where $G$ is metrizable and $\alpha$ is a $d$-open
action) we have.

\begin{cor}\label{cor4.2}
Let $\chi(X)\leq\aleph_0$.
\begin{itemize}
\item[(a)] If $X$ is a coset space of a subgroup of a product of \v Cech complete groups (or there is a transitive and $d$-open action on $X$ of
a subgroup of a product of \v Cech complete groups), then $X$ with
the maximal equiuniformity is an inframetrizable space~\cite[Remark
2.6]{K2013} as a coset space of an inframetrizable group (or since
$X$ admits a transitive and $d$-open action of an inframetrizable
group).
\item[(b)] If $X$ is a coset space of an $\omega$-balanced group (or there is a transitive and $d$-open action on $X$ of
an $\omega$-balanced group), then $X$ is metrizable as a coset space
of a metrizable group (or since $X$ admits a $d$-open action of a
metrizable group).
\item[(c)] If $X$ is a coset space of an $\omega$-narrow group (or there is a transitive and $d$-open action on $X$ of
an $\omega$-narrow group), then $X$ is separable metrizable as a
coset space of a separable metrizable group (or since $X$ admits a
transitive $d$-open action of a separable metrizable group).
\end{itemize}
\end{cor}

Applying theorem of extension of action to the completion of a phase
space with respect to an equiunuformity~\cite[Theorem 3.1]{Mg} and
in the case of actions of $\omega$-narrow groups compactification
theorem~\cite[Theorem 2.13]{Mg2} we have.

\begin{cor}\label{cor4.21}
Let $\chi(X)\leq\aleph_0$. If $X$ is a coset space of an
$\omega$-balanced (respectively $\omega$-narrow) group $G$ or there
is a continuous transitive $d$-open action of an $\omega$-balanced
(respectively $\omega$-narrow) group $G$, then there exists an
equiuniformity (respectively totally bounded equiuniformity)
$\mathcal U$ on $X$ such that $w({\mathcal U})\leq\aleph_0$.

Moreover, there exists an equivariant extension $(G, \tilde X,
\alpha)$ of $(G, X, \alpha)$ and a complete metric on $\tilde X$ for
which the induced uniformity is an extension of $\mathcal U$
(respectively $G$-compactification $(G, bX, \tilde\alpha)$ of $(G,
X, \alpha)$ where $bX$ is a  metrizable compactum).
\end{cor}

\begin{ex}\label{ex4.2}
{\rm Since the Sorgenfrey line is not metrizable, there is no
$d$-open transitive action of an $\omega$-balanced group on the
Sorgenfrey line. Hence, the Sorgenfrey line is not a coset space of
an $\omega$-balanced group.

The Sorgenfrey line is not \v Cech complete. Since a $d$-open action
of a \v Cech complete group is open~\cite[Theorem 3]{K4} and a coset
space of a \v Cech complete group is \v Cech complete~\cite[Theorem
2]{Brown}, there is no $d$-open transitive action of a \v Cech
complete group on the Sorgenfrey line.}
\end{ex}

\begin{que}
Can the Sorgenfrey line be a coset space of a subgroup of a product
of \v Cech complete groups (or, equivalently, of an inframetrizable
group)?
\end{que}

\begin{thm}\label{pro4.1}
Let $X$ be a compact coset space of a subgroup of a product of \v
Cech complete groups (or there is a continuous transitive $d$-open
action of a subgroup of a product of \v Cech complete groups),
$\chi(X)\leq\aleph_0$. Then $X$ is metrizable.
\end{thm}

{\bf Proof.}\quad By item (a) of corollary~\ref{cor4.2} $X$ is a
coset space of an inframetrizable group $G$ (or there is a
continuous transitive $d$-open action of an inframetrizable group on
$X$). Then $(G, X, \alpha)$ is a $G$-Tychonoff space with a natural
open action $\alpha$ in the first case and $d$-open action in the
latter case. By~\cite[Theorem 3.1]{Mg} the extension
$\tilde\alpha:\rho G\times X\to X$ of action $\alpha$ is well
defined and $d$-open. The Raikov completion $\rho G$ is a \v Cech
complete group and by~\cite[Theorem 3]{K4} the action $\tilde\alpha$
is open. Hence, $X$ is a coset space of a \v Cech complete group.
By~\cite[Propositions 3.2.1, 3.2.2]{Chob} $X$ is also a coset space
of an $\omega$-narrow group. Item (c) of corollary~\ref{cor4.2}
finishes the proof. \ $\Box$

\begin{ex}\label{ex4.4}
The two-arrows space~\cite[Exercizes 3.10.C]{E} is not a coset space
of a subgroup of a product of \v Cech complete groups.
\end{ex}

\begin{rem}
{\rm From the proof of theorem~\ref{pro3.1} we can deduce a positive
answer on~\cite[Question 5.7]{K2013} in case of compacta $X$ with
$\chi (X)\leq\aleph_0$.}
\end{rem}

\begin{cor}\label{cor4.22}
Let $\chi(X)\leq\aleph_0$. There is no continuous, effective and
transitive action on $X$ of an $\omega$-balanced group $G$ with
$\psi(G)>\aleph_0$.
\end{cor}

\begin{cor}\label{cor4.11}
Let $X$ be a homogeneous space with Baire property,
$\chi(X)\leq\aleph_0$. If there is a continuous transitive action of
an $\omega$-narrow group on $X$, then $X$ is a separable metrizable
space.
\end{cor}

{\bf Proof.}\quad In~\cite{Usp} it is proved that a transitive
action $\alpha$ of an $\omega$-narrow group $G$ on a space $X$ with
Baire property is $d$-open. Hence, $(G, X, \alpha)$ is a
$G$-Tychonoff space with a transitive $d$-open action. By item (c)
of corollary~\ref{cor4.2} $X$ admits a transitive $d$-open action of
a separable metrizable group. In~\cite[Corollary 4]{CK3} it is shown
that a $d$-open action of a metrizable group yields that the space
$X$ is metrizable. Moreover, if the acting group $G$ is separable
and the action is transitive then the phase space $X$ is separable.\
$\Box$

\begin{ex}\label{ex4.3}
{\rm There is no continuous and transitive action of an
$\omega$-narrow group on the two-arrows space and the homogeneous
compactum $X$ of V.\,V.~Fedorchuk~\cite{Fed} which is not a coset
space and $\chi(X)=\aleph_0$ (see, also, \cite{Chat}).

The first result can also be deduced from~\cite{Usp}, since the
two-arrows space is not dyadic.

\medskip

The Sorgenfrey line has the Baire property. There is no continuous
and transitive action of an $\omega$-narrow group on the Sorgenfrey
line.}
\end{ex}

\begin{cor}\label{cor4.14}
Let $(G, X, \alpha)$ be a $G$-Tychonoff space and $\mathcal U$ be an
equiuniformity on $X$ such that $w({\mathcal U})\leq\aleph_0$. Then
there exist: a $G$-Tychonoff space $(H, X, \gamma)$ such that
\begin{itemize}
\item[(a)] $H$ is an inframetrizable group if $G$ is a subgroup of a product of \v Cech complete
groups,
\item[(b)] $H$ is a metrizable group if $G$ is an $\omega$-balanced
group,
\item[(c)] $H$ is a separable metrizable group if $G$ is an $\omega$-narrow
group;
\end{itemize}
$\mathcal U$ is an equiuniformity on $X$ in $(H, X, \gamma)$ and

\noindent an equivariant pair of maps $(\varphi, {\rm id}):(G, X,
\alpha)\to (H, X, \gamma)$, where $\varphi$ is an epimorphism.
\end{cor}

Applying theorem of extension of action to the completion of a phase
space with respect to an equiunuformity~\cite[Theorem 3.1]{Mg} we
have.

\begin{cor}\label{cor4.15}
Let $(G, X, \alpha)$ be a $G$-Tychonoff space and $\mathcal U$ be an
equiuniformity on $X$ such that $w({\mathcal U})\leq\aleph_0$. Then
there exist: an action $\gamma: H\times\tilde X\to\tilde X$, where
$\tilde X$ is the (metrizable) completion of $X$ with respect to
$\mathcal U$;
\begin{itemize}
\item[(a)]  $H$ is a \v Cech complete group if $G$ is a subgroup of a product of \v Cech complete
groups,
\item[(b)]  $H$ is a complete metrizable group if $G$ is an $\omega$-balanced
group,
\item[(c)]  $H$ is a Polish group if $G$ is an $\omega$-narrow
group,
\end{itemize}
such that $(H, \tilde X, \gamma)$ is a $G$-Tychonoff space, the
extension $\tilde{\mathcal U}$ of $\mathcal U$ on $\tilde X$ is an
equiuniformity in $(H, \tilde X, \gamma)$; and

\noindent an equivarinat embedding $(\varphi, {\rm in}):(G, X,
\alpha)\to (H, \tilde X, \gamma)$, where $\varphi$ is an epimorphism
and ${\rm in}$ is a natural embedding.
\end{cor}

\begin{rem}\label{rem4.16}
{\rm If in corollary~\ref{cor4.15} uniformity $\mathcal U$ is
complete, then $\tilde X=X$.}
\end{rem}

\begin{cor}
Let $(G, X, \alpha)$ be a $G$-Tychonoff space, $G$ be an
$\omega$-narrow group and $\mathcal U$ be an equiuniformity on $X$
such that $w({\mathcal U})\leq\aleph_0$. Then there exist a totally
bounded equiuniformity on $X$ such that $w({\mathcal
U})\leq\aleph_0$ and, hence, an equivariant compactification $(G,
bX, \alpha)$ of $(G, X, \alpha)$, where $bX$ is a metrizable
compactum.
\end{cor}

\begin{ex}\label{ex3.7}
{\rm In~\cite{vM} J.~van Mill constructed a homogeneous Polish space
$Z$ which is not a coset space and proved that no $\omega$-narrow
group acts transitively on $Z$ by a separately continuous
action~\cite[Corollary 3]{vM}.

A) There is no complete metric on $Z$ for which the group $G$ of
uniform equivalences in the topology of uniform convergence with
respect to which acts $d$-openly and has one component of action. In
fact, in this case the group $G$ is Raikov complete~\cite[Chapter X,
\S\ 3, Problem 16 c)]{Burb} and, hence, is complete metrizable. Its
action is open by~\cite[theorem 3]{K4} and has one component of
action. Hence, $Z$ is a coset space but this is a contradiction.

B) There is no totally bounded metric the group $G$ of uniform
equivalences in the topology of uniform convergence with respect to
which acts transitively.  In fact, in this case $G$ is separable
metrizable~\cite{AAK}. This is a contradiction with the properties
of $Z$.

C) Let a subgroup $G$ of a product of \v Cech complete groups acts
continuously, $d$-openly and with one component of action on $Z$.
Then there is no complete equiuniformity with a countable base on
$Z$. In fact, then by corollary~\ref{cor4.15} and
remark~\ref{rem4.16} there is a $d$-open action on $Z$ with one
component of action of a \v Cech complete group. Its action is open
by~\cite[theorem 3]{K4} and has one component of action. Hence, $Z$
is a coset space but this is a contradiction.}
\end{ex}


\end{document}